\documentstyle[twoside,amstex]{article}
\input amssym.def
\input amssym.tex
\def\bc{\begin{center}}
\def\ec{\end{center}}
\def\no{\noindent}
\def\hang{\hangindent\parindent}
\def\textindent#1{\indent\llap{[#1]\enspace}\ignorespaces}
\def\re{\par\hang\textindent}
\begin{document}
\thispagestyle{empty} \vspace*{3 true cm} \pagestyle{myheadings}
\markboth {\hfill {\sl Huanyin Chen and Marjan Sheibani}\hfill}
{\hfill{\sl Feckly Adequate Conditions and Elementary Matrix Reduction}\hfill} \vspace*{-1.5 true cm} \bc{\large\bf Feckly Adequate Conditions and Elementary Matrix Reduction}\ec

\vskip6mm
\bc{{\bf Huanyin Chen}\\[2mm]
Department of Mathematics, Hangzhou Normal University\\
Hangzhou 310036, China\\
huanyinchen@@aliyun.com}\ec

\bc{{\bf Marjan
Sheibani}\\[2mm]
Faculty of Mathematics, Statistics and Computer Science\\
Semnan University, Semnan, Iran\\
m.sheibani1@@gmail.com}\ec

\vskip10mm
\begin{minipage}{120mm}
\no {\bf Abstract:} We present some new conditions for a B$\acute{e}$zout ring to be an elementary divisor ring.
We prove, in this note, that a B$\acute{e}$zout ring $R$ is feckly zero-adequate if and only if
$R/J(R)$ is regular if and only if $R/J(R)$ is $\pi$-regular, and that every feckly zero-adequate ring is an elementary divisor ring. If $R$ has feckly adequate range 1, we prove that $R$ is an elementary divisor ring if and only if $R$ is a B$\acute{e}$zout ring. Many known results are
thereby generalized to much wider class of rings, e.g. [4, Theorem 14], [5, Theorem 4], [8, Theorem 1.2.14], [10, Theorem 4] and [11, Theorem 7].
\vskip3mm \no {\bf MR(2010) Subject Classification}: 13F99, 13E15,
06F20.
\end{minipage}

\vskip15mm \bc{\bf 1. Introduction}\ec

\vskip4mm \no Throughout this paper, all rings are commutative
with an identity. A matrix $A$ (not necessarily square) over a ring
$R$ admits diagonal reduction if there exist invertible matrices
$P$ and $Q$ such that $PAQ$ is a diagonal matrix $(d_{ij})$, for
which $d_{ii}$ is a divisor of $d_{(i+1)(i+1)}$ for each $i$. A
ring $R$ is called an elementary divisor ring provided that every
matrix over $R$ admits a diagonal reduction. A ring $R$ is a Hermite ring if every $1\times 2$ matrix over $R$
admits a diagonal reduction. As is well known, a ring $R$ is
Hermite if and only if for all $a,b\in R$ there exist $a_1,b_1\in
R$ such that $a=a_1d,b=b_1d$ and $a_1R+b_1R=R$ ([8, Theorem 1.2.5]).
A ring is a
B$\acute{e}$zout ring if every finitely generated ideal is
principal. Obviously, $\{ ~\mbox{elementary divisor rings}~\}\subsetneq \{ ~~\mbox{Hermite rings}~\}\subsetneq \{ ~~\mbox{B$\acute{e}$zout rings}~\} ( ~\mbox{cf. [8]}).$
An attractive problem is to investigate various conditions under which a B$\acute{e}$zout ring is an elementary divisor ring.

We recall that an element $c\in R$ is adequate provided that for any $a\in R$ there exist some $r,s\in R$ such that $(1)$ $c=rs$; $(2)$ $rR+aR=R$; $(3)$ $s'R+cR\neq R$ for each non-invertible divisor $s'$ of $s$. Whether a ring with various adequate properties is an elementary ring is studied by many authors.
A ring $R$ is clean provided that every element in $R$ is the sum and a unit. In [10, Theorem 4], Zabavsky and Bilavska proved that an interesting result: every zero-adequate ring, i.e., a B$\acute{e}$zout ring in which zero is adequate, is a clean ring. Recently, Pihua claimed that every zero-adequate ring is semiregular [5, Theorem 4], and so such kind ring is an elementary divisor ring [8, Theorem 2.5.2]. A B$\acute{e}$zout ring is called an adequate ring provided that every nonzero element is adequate. In his research of elementary divisor domains, Helmer proved that every adequate domain is an elementary divisor ring. After his work, Kaplansky showed that an adequate ring whose zero divisors are in the radical is an elementary divisor ring. Helmer also showed that an adequate ring is an elementary divisor ring if and only it is a Hermite ring. For more results about adequate conditions, we refer the reader to [8].

Recall that a ring $R$ has stable range 1 if $aR+bR=R$ with $a,b\in R$ there exists a $y\in R$ such that $a+by\in R$ is invertible. Such condition plays an important role in algebraic K-theory (cf. [2]). It includes many kind of rings, e.g., regular rings, semiregular rings, $\pi$-regular rings, local rings, clean rings, etc. Domsha and Vasiunyk combined this condition with adequate condition together. A ring $R$ is called to have adequate range 1 if $aR+bR=R$ with $a,b\in R$ implies that there exists a $y\in R$ such that $a+by\in R$ is adequate. It was proved that every B$\acute{e}$zout domain having adequate range 1 is an elementary divisor ring [4, Theorem 14].

In this note, we are concern on a new condition related to the Jacobson radical of a ring.
We say that $c\in R$ is feckly adequate if for any $a\in R$ there exist
some $r,s\in R$ such that $(1)$ $c\equiv rs~(mod~J(R))$; $(2)$ $rR+aR=R$; $(3)$
$s'R+aR\neq R$ for each non-invertible divisor $s'$ of $s$. A B$\acute{e}$zout ring $R$ is called a feckly zero-adequate ring provided that $0\in R$ is feckly adequate.
In Section 2, we investigate the necessary and sufficient conditions under which a ring $R$ is feckly zero-adequate. We prove,
a B$\acute{e}$zout ring $R$ is feckly zero-adequate if and only if
$R/J(R)$ is regular if and only if $R/J(R)$ is $\pi$-regular; hence that $R$ is an elementary divisor ring. Examples of feckly zero-adequate rings which are not zero-adequate are provided. Further, we shall characterize feckly zero-adequate rings in terms of the
idempotent stable range $1$. In Section 3, we explore a new condition behaving like the stable range under which a B$\acute{e}$zout ring is an elementary divisor ring. A ring $R$ is said to have feckly adequate range 1 provided that $aR+bR=R$ implies that there exists a $y\in R$ such that $a+by\in R$ is feckly adequate. If $R$ has feckly adequate range 1, we prove that $R$ is an elementary divisor ring if and only if $R$ is a B$\acute{e}$zout ring. This extend [8, Theorem 1.2.14] and [4, Theorem 14] to much wider
class of rings.

We shall use $J(R)$ and $U(R)$ to denote the Jacobson radical of $R$ and the set of all units in $R$, respectively. A ring $R$ is called a domain if there is no any nonzero zero divisor of $R$.

\vskip15mm\bc{\bf 2. Feckly Zero-adequate Rings}\ec

\vskip4mm The purpose of this section is to characterize feckly zero-adequate rings. A ring $R$ is $\pi$-regular if for any $a\in R$ there exists $n\in {\Bbb N}$ such that $a^n=a^nba^n$ for some $b\in R$. We begin with

\vskip4mm \hspace{-1.8em} {\bf Lemma 2.1.}\ \ {\it Let $R$ be a ring. Then $R/J(R)$ is $\pi$-regular if and only if for any $a\in R$, there exists an $n\in {\Bbb N}$, an element $e\in R$, a unit $u\in R$ and a $w\in J(R)$ such that $a^n=eu+w$ and $e-e^2\in J(R)$.}
\vskip2mm\hspace{-1.8em} {\it Proof.}\ \ $\Longrightarrow$ Let $a\in R$. Then $\overline{a^n}=\overline{a^nba^n}$ for some $n\in {\Bbb N}$.
Set $e=a^nb$ and $u=1-a^nb+a^n$. Then $\overline{e^2}=\overline{e}\in R/J(R)$ and $\big(\overline{u}\big)^{-1}=1-a^nb+ba^nb$ in $R/J(R)$. As units lift modulo $J(R)$,
 we see that $u\in U(R)$. Set $w:=a-eu$. Then we obtain $a=eu+w$, where $e^2-e,w\in J(R)$.

$\Longleftarrow$ For any $a\in R$, there exists an $n\in {\Bbb N}$, an element $e\in R$, a unit $u\in R$ and a $w\in J(R)$ such that $a^n=eu+w$ and $e-e^2\in J(R)$. Hence, $\overline{a^n}=\overline{eu}$ in $R/J(R)$. Therefore $\overline{a^n}=\overline{a^nu^{-1}a^n}$, as required.\hfill$\Box$

\vskip4mm \hspace{-1.8em} {\bf Lemma 2.2.}\ \ {\it Let $R$ be a B$\acute{e}$zout ring. If $R/J(R)$ is $\pi$-regular, then $R$ is feckly zero-adequate.}
\vskip2mm\hspace{-1.8em} {\it Proof.}\ \ Suppose that $R/J(R)$ is $\pi$-regular.
Let $a\in R$ be an arbitrary element. In light of Lemma 2.1, there exists an element $e\in R$, a unit $u\in R$ and a $w\in J(R)$ such that
$a^n=eu+w (n\in {\Bbb N})$ and $e-e^2\in J(R)$. Then $(1-e)e\in J(R)$. Clearly, $a^nu^{-1}+(1-e)=1+wu^{-1}\in U(R)$, and so $a^nu^{-1}\big(1+wu^{-1}\big)^{-1}+(1-e)\big(1+wu^{-1}\big)^{-1}=1$.
Hence, $(1-e)R+a^nR=R$. This implies that $(1-e)R+aR=R$. If $s$ is a non-invertible divisor of $e$, then $e=ss'$ for some $s'\in R$.
If $sR+aR=R$, then $sR+a^nR=R$. Thus, we can find some $x,y\in R$ such that
$sx+a^ny=1$. Hence, $sx+(eu+w)y=1$, and so $s(x+s'uy)=1-wy\in U(R)$. This implies that $s$ is invertible, an absurd.
Therefore $sR+aR\neq R$. Accordingly, $R$ is feckly zero-adequate.\hfill$\Box$

\vskip4mm Recall that a ring is feckly clean provided that for any $a\in R$ there exists an element $e\in R$ such that $a-e\in U(R)$ and $e-e^2\in J(R)$.
As is well known, a ring $R$ is feckly clean if and only if $aR+bR=R$ with $a,b\in R$ implies that there are $x,y\in R$ such that $a|x, b|y, xy\in J(R)$ and $xR+yR=R$, if and only if $Max(R)$ is zero-dimensional ([6, Theorem 3.13 and Proposition 3.12]). More topological characterizations of such type of rings, we refer the reader to [3].

\vskip4mm \hspace{-1.8em} {\bf Lemma 2.3.}\ \ {\it Every feckly zero-adequate ring is feckly clean.}
\vskip2mm\hspace{-1.8em} {\it Proof.}\ \ Let $R$ be a feckly zero-adequate ring. Let $x\in
R$. Then we have some $r,s\in R$ such that $rs\in J(R)$,
where $rR+xR=R$ and $s'R+xR\neq R$ for any
noninvertible divisor $s'$ of $s$. We claim that
$rR+sR=R$. If not, $rR+sR=hR\neq R$. Thus, $h\in R$ is
a noninvertible divisor of $s$; hence that $hR+xR\neq
R$. But $h$ is a divisor of $r$, we get
$hR+xR=R$, which is impossible. Write
$rc+sd=1$ in $R$. Then $(rc)^2-rc=(rc)^2-rc\big(rc+sd)=-(rs)(cd)\in J(R)$. Set
$e=rc$. Then $e^2-e\in J(R)$.

Claim I. $(x-e)R+rR=R$. If not, $(x-e)R+rR=tR\neq R$. Then
$rR\subseteq tR$, and so $xR\subseteq tR$. This implies that
$rS+xS\subseteq tR\neq R$, a contradiction.

Claim II. $(x-e)R+sR=R$. If not, $(x-e)R+sR=tR\neq R$. Then $t$ is
a noninvertible divisor of $s$, and so $tR+xR\neq R$.
But $eR+sR=R$, and so $eR+tR=R$. Write $x-e=tw$ with $w\in R$.
Then $e=x-tw$, and so $eR+tR\subseteq tR+xR\neq R$, which is
impossible.

Therefore $(x-e)R+rsR=R$. Write $(x-e)p+(rs)q=1$ for some $p,q\in R$. As $rs\in J(R)$, we deduce that $(x-e)p=1-(rs)q\in U(R)$, and so $x-e\in U(R)$. This completes the proof.\hfill$\Box$

\vskip4mm \hspace{-1.8em} {\bf Lemma 2.4.}\ \ {\it Let $R$ be a feckly zero-adequate ring. Then $J(R)=\{ x~|~x-u\in U(R)~\mbox{for any}$ $u\in U(R)\}$.}
\vskip2mm\hspace{-1.8em} {\it Proof.}\ \ Clearly, $J(R)\subseteq \{ x~|~x-u\in U(R)~\mbox{for any}~u\in U(R)\}$. Let $x\in R$ and $x-u\in U(R)$ for any $u\in U(R)$. Let $r\in R$. Then $xR+(1-xr)R=R$. Since $R$ is feckly clean, by Lemma 2.3, $R/J(R)$ is clean, and then $R/J(R)$ has stable range 1. It follows that
 $R$ has stable range 1. Thus, we have a $y\in R$ such that $u:=x+(1-xr)y\in U(R)$. Hence, $x-u=-(1-xr)y\in U(R)$, and then $1-xr\in U(R)$. Therefore $x\in J(R)$, as desired.\hfill$\Box$

\vskip4mm Recall that a ring $R$ is regular if for any $a\in R$ there exists a $b\in R$ such that $a=aba$. We are now ready to prove:

\vskip4mm \hspace{-1.8em} {\bf Theorem 2.5.}\ \ {\it Let $R$ be a B$\acute{e}$zout ring. Then the following are equivalent:}
\begin{enumerate}
\item [(1)] {\it $R$ is feckly zero-adequate.}
\vspace{-.5mm}
\item [(2)] {\it $R/J(R)$ is regular.}
\vspace{-.5mm}
\item [(3)] {\it $R/J(R)$ is $\pi$-regular.}
\end{enumerate}
\vspace{-.5mm} {\it Proof.}\ \ $(1)\Rightarrow (2)$ Let $x\in
R$. Then we have some $r,s\in R$ such that $rs\in J(R)$,
where $rR+xR=R$ and $s'R+xR\neq R$ for any
noninvertible divisor $s'$ of $s$. As in the proof of Lemma 2.3, we see that
$rR+sR=R$. Since $rR+xR=rR+sR=R$, we get $rR+sxR=R$.
Write
$rc+sxd=1$ in $R$. Set
$e=rc$. Then $e^2-e=(rc)^2-rc=(rc)(sxd)\in J(R)$. Let $u$ be an arbitrary invertible element of $R$.

Claim I. $(u-ex)R+rR=R$. If not, $(u-ex)R+rR=tR\neq R$. Then
$rR\subseteq tR$, and so $u\in eR+tR\subseteq rR+tR\subseteq tR$. This implies that
$t\in U(R)$, a contradiction.

Claim II. $(u-ex)R+sR=R$. If not, $(u-ex)R+sR=tR\neq R$. Then $t$ is
a noninvertible divisor of $s$, and so $tR+xR\neq R$.
It follows from $eR+sR=R$ that $eR+tR=R$. Write $u-ex=tw$ with $w\in R$.
Then $e=ex+tw$, and so $eR+tR\subseteq tR+xR\neq R$, an absurd.

Finally, $(u-ex)R+rsR=R$. As $rs\in J(R)$, we get $u-ex\in U(R)$. This implies that
$\big(x-(1-e)x\big)-u=ex-u\in U(R)$. In view of Lemma 2.3, $R$ has stable range 1.
It follows by Lemma 2.4 that $x-(1-e)x\in J(R)$. Clearly, $1-e=sxd\in xR$. Therefore $\overline{x}=\overline{x(sd)x}$ in $R/J(R)$, and so $R/J(R)$ is regular.

$(2)\Rightarrow (3)$ This is obvious.

$(3)\Rightarrow (1)$ In view of Lemma 2.2, $R$ is feckly zero-adequate, as asserted.\hfill$\Box$

\vskip4mm \hspace{-1.8em} {\bf Corollary 2.6.}\ \ {\it Every feckly zero-adequate ring is an elementary divisor ring.}
\vskip2mm\hspace{-1.8em} {\it Proof.}\ \ Let $R$ be a feckly zero-adequate ring. Then $R/J(R)$ is regular, by Theorem 2.5. In view of [8, Theorem 1.2.14 and Theorem 1.2.13], $R/J(R)$ is an elementary divisor ring. Therefore $R$ is an elementary divisor ring.\hfill$\Box$

\vskip4mm \hspace{-1.8em} {\bf Corollary 2.7.}\ \ {\it A ring $R$ is feckly zero-adequate if and only if}
\begin{enumerate}
\item [(1)] {\it $R$ is a B$\acute{e}$zout ring;}
\vspace{-.5mm}
 \item [(2)] {\it $R/J(R)$ is zero-adequate.}
\end{enumerate}
\vspace{-.5mm} {\it Proof.}\ \ $\Longrightarrow $ $(1)$ is obvious. In view of Theorem 2.5, $R/J(R)$ is regular, proving $(2)$, as every element in $R/J(R)$ is adequate.

$\Longleftarrow $ For any $a\in R$, we can find some $r,s\in R$ such that $\overline{0}=\overline{rs}$,
where $\overline{r}\big(R/J(R)\big)+\overline{x}\big(R/J(R)\big)=R/J(R)$ and $\overline{s'}\big(R/J(R)\big)+\overline{x}\big(R/J(R)\big)\neq R/J(R)$ for any
noninvertible divisor $s'$ of $s$. Thus, $rs\in J(R)$. Further,
$rR+xR=R$ and $s'R+xR\neq R$ for any
noninvertible divisor $s'$ of $s$. Therefore $R$ is feckly zero-adequate.\hfill$\Box$

\vskip4mm We note that $"\Longrightarrow "$ in Corollary 2.7 can not be proved in a direct route by the definitions. Let $R={\Bbb Z}[\alpha]$, where $\alpha^2=1$. Choose $J=(1+\alpha), s'=5-3\alpha,s=3+\alpha\in R$. Then $\overline{s'}=2$ is a noninvertible divisor of $\overline{s}=4$ in $R/J$, while $s'$ is not a noninvertible divisor of $s$ in $R$.

\vskip4mm \hspace{-1.8em} {\bf Corollary 2.8.}\ \ {\it Let $R$ be a B$\acute{e}$zout ring. Then $R$ is zero-adequate if and only if}
\begin{enumerate}
\item [(1)] {\it $R$ is feckly zero-adequate;}
\vspace{-.5mm}
\item [(2)] {\it Every idempotent lifts modulo $J(R)$.}
\end{enumerate}
\vspace{-.5mm} {\it Proof.}\ \ $\Longrightarrow $ $(1)$ is obvious. In view of [10, Theorem 4], $R$ is clean, proving $(2)$, by the Nicolson Theorem.

$\Longleftarrow $ In view of Theorem 2.5, $R/J(R)$ is regular. As idempotent lifts modulo $J(R)$, as in the proof of Lemma 2.2, $R$ is zero-adequate.\hfill$\Box$

\vskip4mm A ring $R$ is semiregular if $R/J(R)$ is regular and idempotents lift modulo $J(R)$. We now derive

\vskip4mm \hspace{-1.8em} {\bf Corollary 2.9 [5, Theorem 4].}\ \ {\it Let $R$ be a B$\acute{e}$zout ring. Then $R$ is zero-adequate if and only if $R$ is semiregular.}
\vskip2mm\hspace{-1.8em} {\it Proof.}\ \ This is obvious, by Corollary 2.8 and Theorem 2.5.\hfill$\Box$

\vskip4mm From Corollary 2.9, we observe that the difference between feckly adequate rings and adequate rings is just to forget the lifting of idempotents.

\vskip4mm \hspace{-1.8em} {\bf Example 2.10.}\ \ {\it Every finite B$\acute{e}$zout ring is feckly zero-adequate.}
\vskip2mm\hspace{-1.8em} {\it Proof.}\ \ Since every finite ring is $\pi$-regular, we complete the proof, by Theorem 2.5.\hfill$\Box$

\vskip4mm \hspace{-1.8em} {\bf Example 2.11.}\ \ {\it Let $R=\{ \frac{m}{n}~\mid~m,n\in {\Bbb Z}, n\neq 0, 3,5 \nmid n\}$. Then $R$ is feckly zero-adequate, while $R$ is not zero-adequate.}
\vskip2mm\hspace{-1.8em} {\it Proof.}\ \ Let $I=\frac{a}{b}R+\frac{c}{d}R$, where $\frac{a}{b},\frac{c}{d}\in R$. As ${\Bbb Z}$ is a principal ideal domain, we can find some $p\in {\Bbb Z}$ such that $a{\Bbb Z}+c{\Bbb Z}=p{\Bbb Z}$. One easily checks that $I=pR$. Thus, $R$ is a B$\acute{e}$zout ring. As in the proof of [1, Example 17], $R$ has only
two maximal ideals $3R$ and $5R$. Since $3R+5R=R$, by Chinese Reminder Theorem, we deduce that $R/J(R)\cong R/3R \times R/5R\cong {\Bbb Z}_3\times {\Bbb Z}_5$.
As ${\Bbb Z}_3$ and ${\Bbb Z}_5$ are regular rings, $R/J(R)$ is regular. It follows by Theorem 2.5 that $R$ is feckly zero-adequate. As in the proof of [1, Example 17], $R$ is not clean. In light of [10, Theorem 4], $R$ is not zero-adequate, as desired.\hfill$\Box$

\vskip4mm \hspace{-1.8em} {\bf Example 2.12.}\ \ {\it Let $F$ be
a field, and let $R=F[[x,y]]$. Let $S=R-
(x)\bigcup (y)$. Then $R_S$ is feckly zero-adequate, but it is
not zero-adequate.} \vskip2mm\hspace{-1.5em}  {\it Proof.}\ \ As
$F[[x,y]]/(x)\cong F[[y]]$ is an integral domain, we see that
$(x)$ is a prime ideal of $F[[x,y]]$. Likewise, $(y)$ is a prime
ideal of $R$. Set $S=R-(x)\bigcup (y)$. Then $S$ is a
multiplicative closed subset of $R$. Let $P$ be a maximal ideal of
$R_S$. Then we can find an ideal $Q$ of $R$ such that $P=Q_S$ such
that $Q\bigcap S=\emptyset$. Hence, $Q\subseteq (x)\bigcup (y)$.
Assume that $Q\nsubseteq (x)$ and $Q\nsubseteq (y)$. Then we can
find some $b\in Q$, but $b\not\in (x)$. Likewise, we have some
$c\in Q$, but $c\not\in (y)$. Set $a=b+c$. Then $a\in Q$, but
$a\not\in (x)\bigcup (y)$. This gives a contradiction. Hence,
$Q\subseteq (x)$ or $Q\subseteq (y)$. It follows that $Q_S\subseteq
(x)_S$ or $Q_S\subseteq (y)_S$. By the maximality of $P$, we get
$P=(x)_S$ or $(y)_S$. Thus, $R_S$ has exactly two maximal
ideals $(x)_S$ and $(y)_S$. Accordingly, $R_S/J(R_S)\cong R_S/(x)_S\times R_S/(y)_S$ is regular. Obviously, $R_S$ is a PID, and then it is a B$\acute{e}$zout ring. In light of Theorem 2.5,
$R_S$ is feckly zero-adequate. Obviously, $R_S$ is indecomposable. This implies that $R_S$ is not clean; otherwise,
it is local, which is impossible. Therefore $R_S$ is not zero-adequate, by [10, Theorem 4].\hfill$\Box$

\vskip4mm A B$\acute{e}$zout ring in which every nonzero element is adequate
is called an adequate ring. It is an open problem if a homomorphic image of an adequate ring is again adequate (cf. [9, Question 32]). We now give an affirmative answer for a similar problem on feckly zero-adequate rings.

\vskip4mm \hspace{-1.8em} {\bf Proposition 2.13.}\ \ {\it Every homomorphic image of a feckly zero-adequate ring is feckly zero-adequate.}
\vskip2mm\hspace{-1.8em} {\it Proof.}\ \ Let $I$ be an ideal of feckly zero-adequate ring $R$. We shall prove that $R/I$ is feckly zero-adequate. As every homomorphic image of a B$\acute{e}$zout ring is B$\acute{e}$zout, then by Therem 2.5 we need only to prove that $\big(R/I\big)/J\big(R/I\big)$ is a regular ring. Let $\overline{a} \in R/I$. Since $R$ is feckly zero-adequate, then by Theorem 2.5, $R/J(R)$ is regular.
In light of Lemma 2.1. there exists an element $e$, a unit $u\in R$ and a $w\in J(R)$ such that $a=eu+w$ and $e-e^2\in J(R)$. Now we have $\overline{a}=\overline{e}\overline{u}+\overline{w}$ in $R/I$. It is obvious that $\overline{e}^2-\overline{e}\in J(R/I)$, $\overline{u}$ is a unit of $R/I$ and $\overline{w}\in J(R/I)$. In virtue of Lemma 2.1, $\big(R/I\big)/J\big(R/I\big)$ is regular. Therefore $R/I$ is feckly zero-adequate, by Theorem 2.5.\hfill$\Box$

\vskip4mm \hspace{-1.8em} {\bf Corollary 2.14.}\ \ {\it Let $R=\prod\limits_{i\in I}R_{i}$. Then $R$ is feckly zero-adequate if and only if each
$R_i~(i\in I)$ is feckly zero-adequate.}
\vskip2mm\hspace{-1.8em} {\it Proof.}\ \ Let $R$  be a feckly zero-adequate ring. By Proposition 2.13, every homomrphic image of $R$ is feckly zero-adequate,
and then $R_i$ is feckly zero-adequate for each $i\in I$. Now let each $R_i~(i\in I)$ be feckly zero-adequate. Then each $R_i$ is a B$\acute{e}$zout ring, and thus, $R$ is a B$\acute{e}$zout ring. In view of Theorem 2.5, $R_i/J(R_i)$ is regular. For any $(a_i)\in R$, we have $a_i\in R_i$, and so there exists an element $b_i\in R_i$ such that $a_i-a_ib_ia_i\in J(R_i)$. Hence, $(a_i)-(a_i)(b_i)(a_i)\in J(R)$. This implies that $R/J(R)$ is regular. Accordingly, $R$ is feckly zero-adequate, by Theorem 2.5.
\hfill$\Box$

\vskip4mm \hspace{-1.8em} {\bf Corollary 2.15.}\ \ {\it Let $R$ be a B$\acute{e}$zout ring. Then $R$ is a feckly zero-adequate ring if and only if $R/J(R)$ is feckly zero-adequate.}
\vskip2mm\hspace{-1.8em} {\it Proof.}\ \ $\Longrightarrow $ This is clear by Proposition 2.13.

$\Longleftarrow$ In view of Theorem 2.5, $R/J(R)$ is regular. Since $R$ is a B$\acute{e}$zout ring, by Theorem 2.5 again, $R$ is feckly zero-adequate.\hfill$\Box$

\vskip4mm In [10, Theorem 4], Zabavsky and Bilavska proved that every zero-adequate ring has idempotent stable range 1. We now extend this result and characterize feckly zero-adequate rings as follows.

\vskip4mm \hspace{-1.8em} {\bf Theorem 2.16.}\ \ {\it Let $R$ be a B$\acute{e}$zout ring. Then the following are equivalent:}
\vspace{-.5mm}
\begin{enumerate}
\item [(1)] {\it $R$ is feckly zero-adequate.}
\vspace{-.5mm}
\item [(2)] {\it $aR+bR=R$ with $a,b\in R$ implies that there exists an element $e\in R$ such that $a+be\in U(R), aR\cap eR\subseteq J(R)$ and $e-e^2\in J(R)$.}
\vspace{-.5mm}
\item [(3)] {\it For any $a\in R$, there exists an element $e\in R$ such that $a-e\in U(R), aR\cap eR\subseteq J(R)$ and $e-e^2\in J(R)$.}
\end{enumerate}
\vspace{-.5mm} {\it Proof.}\ \ $(1)\Rightarrow (2)$ Given $aR+bR=R$ with $a,b\in R$, then $\overline{a}\big(R/J(R)\big)+\overline{a}\big(R/J(R)\big)=R/J(R)$.
In view of Theorem 2.5, $R/J(R)$ is regular. As $R$ is commutative, we see that
$R/J(R)$ is unit-regular. By virtue of [2, Corollary 5.3.3], we can find an idempotent $\overline{e}\in R/J(R)$ such that $\overline{a+be}\in U\big(R/J(R)\big)$
and $\overline{a}\big(R/J(R)\big)\cap \overline{e}\big(R/J(R)\big)=\overline{0}$. Set $w=a+be$. As unit lifts modulo $J(R)$, we see that $w\in U(R)$, and then $a+be\in U(R)$. Further, $aR\cap eR\subseteq J(R)$ and $e-e^2\in J(R)$.

$(2)\Rightarrow (3)$ For any $a\in R$, we have $aR+(-1)R=R$. Thus, there exists an element $e\in R$ such that $a-e\in U(R)$ and $aR\cap eR\subseteq J(R), e-e^2\in J(R)$.

$(3)\Rightarrow (1)$ Let $a\in R$. Then there exists an element $e\in R$ such that $u:=a-e\in U(R)$ and $aR\cap eR\subseteq J(R), e-e^2\in J(R)$. It follows that $1=au^{-1}-eu^{-1}$, and so $au^{-1}e=e(1+eu^{-1})\subseteq aR\cap eR$. Hence, $au^{-1}(a-u)=au^{-1}a-a\in J(R)$, and then
$\overline{a}=\overline{au^{-1}a}$.
Thus, $R/J(R)$ is regular. Therefore $R$ is a feckly zero-adequate, by Theorem 2.5.\hfill$\Box$

\vskip4mm \hspace{-1.8em} {\bf Corollary 2.17.}\ \ {\it Let $R$ be a B$\acute{e}$zout ring. Then the following are equivalent:}
\vspace{-.5mm}
\begin{enumerate}
\item [(1)] {\it $R$ is zero-adequate.}
\vspace{-.5mm}
\item [(2)] {\it $aR+bR=R$ with $a,b\in R$ implies that there exists an idempotent $e\in R$ such that $a+be\in U(R), aR\cap eR\subseteq J(R)$.}
\vspace{-.5mm}
\item [(3)] {\it For any $a\in R$, there exists an idempotent $e\in R$ such that $a-e\in U(R), aR\cap eR\subseteq J(R)$.}
\end{enumerate}
\vspace{-.5mm} {\it Proof.}\ \ $(1)\Rightarrow (2)$ In view of Corollary 2.8, $R$ is feckly zero-adequate. Suppose that $aR+bR=R$ with $a,b\in R$. Then
there exists an element $f\in R$ such that $a+bf\in U(R), aR\cap fR\subseteq J(R)$ and $f-f^2\in J(R)$. Further, we can find an idempotent $e\in R$ such that
$e-f\in J(R)$. Hence, $a+be\in U(R)$. If $x\in aR\cap eR$, then $x=ar=es$ for some $r,s\in R$. Then $x=ar=fs+(e-f)s=ex=fes+(e-f)es$. We infer that
$ar-(e-f)es=fes\in aR\cap fR\subseteq J(R)$. It follows that $x=ar\in J(R)$, and therefore $aR\cap eR\subseteq J(R)$.

$(2)\Rightarrow (3)$ This is obvious.

$(3)\Rightarrow (1)$ In view of Theorem 2.16, $R$ is feckly zero-adequate. Furthermore, $R$ is clean, and so every idempotent lifts modulo $J(R)$.
Therefore we complete the proof, by Corollary 2.8.\hfill$\Box$

\vskip15mm \bc{\bf 3. Elementary Matrix Reduction}\ec

\vskip4mm As is well known, an adequate ring is an elementary divisor ring if and only if it is a Hermite ring ([8, Theorem 1.2.14]). The aim of this section is to
extend this result to the rings having feckly adequate range 1.

\vskip4mm \hspace{-1.8em} {\bf Theorem 3.1.}\ \ {\it Let $R$ be a B$\acute{e}$zout ring. If $a\in R$ is feckly adequate, then
$R/aR$ is feckly zero-adequate.} \vskip2mm \hspace{-1.5em}{\it
Proof.}\ \ Let $\overline{b}\in R/aR$. Then there exist $r,s\in R$
such that $a\equiv rs (mod~J(R)), (r,b)=1$ and $(s',b)\neq 1$ for any noninvertible
divisor $s'$ of $s$. Hence, $\overline{a}\equiv \overline{rs} \big(mod~J(R/aR)\big)$, i.e., $\overline{rs}\in J(R/aR)$.
Clearly, $\overline{r}(R/aR)+\overline{b}(R/aR)=R/aR$. Let
$\overline{t}\in R/aR$ be a noninvertible divisor of
$\overline{s}$. Then $t$ is a divisor of $s+ak$ for some $k\in R$. Write
$s+ak=t\beta$ for some $\beta\in R$. Then $s+rsk=t\beta+w$ for a $w\in J(R)$, and so
$s(1+rk)=t\beta+w$. If $sR+tR=R$, then
$sp+tq=1$ for some $p,q\in R$. It follows that
$s(1+rk)p+t(1+rk)q=1+rk$, and so $(t\beta +w)p+t(1+rk)q=1+rk$. As $w\in J(R)$, we get
$r(-k)(1-wp)^{-1}+t(\beta p+(1+rk)q)(1-wp)^{-1}=1$.
This implies
that $rR+tR=R$; hence, $(rs)R+tR=R$. As $a-rs\in J(R)$, we see that $aR+tR=R$, and
then $\overline{t}\in R/aR$ is invertible, a contradiction.
Therefore $sR+tR\neq R$. Since $R$ is a B$\acute{e}$zout ring, we have a noninvertible $u\in R$ such that
$sR+tR=uR$. We infer that $u$ is a noninvertible divisor of $s$.
Hence, $uR+bR\neq R$. This proves that
$\overline{u}(R/aR)+\overline{b}(R/aR)\neq R/aR$; otherwise,
there exist $x,y,z\in R$ such that $ux+by=1+az$. This implies that
$ux+by=1+rsz=1+ucrz$ for a $c\in R$. Hence, $u(x-crz)+by=1$, a
contradiction. Thus $\overline{t}(R/aR)+\overline{b}(R/aR)\neq
R/aR$, and so the result is proved.\hfill$\Box$

\vskip4mm \hspace{-1.8em} {\bf Corollary 3.2.}\ \ {\it Let $R$ be a B$\acute{e}$zout domain. Then $a\in R$ is adequate if and only if}
\vspace{-.5mm}
\begin{enumerate}
\item [(1)] {\it $a\in R$ is feckly adequate;}
\vspace{-.5mm}
\item [(2)] {\it Every idempotent lifts modulo $J(R/aR)$.}
\end{enumerate}
\vspace{-.5mm} {\it Proof.}\ \ $\Longrightarrow $ $(1)$ is obvious. In view of [10, Theorem 2], $\overline{0}\in R/aR$ is adequate. Thus, proving $(2)$ by
Corollary 2.8.

$\Longleftarrow$ In virtue of Theorem 3.1, $R/aR$ is feckly zero-adequate. It follows by Theorem 2.5 that $R/aR/J(R/aR)$ is regular.
Therefore $R/aR$ is semiregular, and so $\overline{0}\in R/aR$ is adequate, by Corollary 2.9. Therefore $a\in R$ is adequate, in terms of
[10, Theorem 10].\hfill$\Box$

\vskip4mm Recall that a B$\acute{e}$zout ring is everywhere adequate provided that every element in $R$ is adequate.

\vskip4mm \hspace{-1.8em} {\bf Proposition 3.3.}\ \ {\it Let $R$ be a B$\acute{e}$zout domain. Then the following are equivalent:}
\vspace{-.5mm}
\begin{enumerate}
\item [(1)] {\it $R$ is zero-adequate;}
\vspace{-.5mm}
\item [(2)] {\it $R$ is everywhere adequate.}
\vspace{-.5mm}
\item [(3)] {\it $R$ is local.}
\end{enumerate}
\vspace{-.5mm} {\it Proof.}\ \ $(1)\Rightarrow (2)$ In view of Corollary 2.9, $R$ is semiregular. Let $a\in R$. Then $R/aR$ is semiregular, by Lemma 2.1.
Hence, $R/aR$ is zero-adequate. It follows by [10, Theorem 10] that $a\in R$ is adequate. Therefore $R$ is everywhere adequate.

$(2)\Rightarrow (3)$ It follows by Corollary 2.9 that $R$ is semiregular. Let $a\in R$.
Then there exists an idempotent $e\in aR$ such that $(1-e)a\in J(R)$. Since
$R$ is a domain, we see that $e=0,1$. If $e=0$, then $a\in J(R)$. If $e=1$, then $a\in U(R)$. Therefore $R$ is local.

$(3)\Rightarrow (1)$ Clearly, $R$ is semiregular, and therefore we complete the proof, by Corollary 2.9.\hfill$\Box$

\vskip4mm \hspace{-1.8em} {\bf Corollary 3.4.}\ \ Every zero-adequate domain is adequate.

\vskip4mm The converse of Corollary 3.4 is not true. For instance, ${\Bbb Z}$ is an adequate domain, but it is not zero-adequate.

\vskip4mm \hspace{-1.8em} {\bf Corollary 3.5.}\ \ {\it Let $R$ be a zero-adequate domain, and let $A=BC$ where $B$ and $C$ are matrices over $R$.
Then elementary divisors of $A$ are divisible by appropriate elementary divisors of $B$ and $C$.}  \vskip2mm \hspace{-1.5em}{\it
Proof.}\ \ In view of Proposition 3.3, $R$ is everywhere adequate. Additionally, $R$ is local, and so it has stable range 1.
Accordingly, we complete the proof by [8, Theorem 4.10.4].\hfill$\Box$

\vskip4mm \hspace{-1.8em} {\bf Lemma 3.6.}\ \
{\it A ring $R$ is an elementary divisor ring if and only if}
\begin{enumerate}
\item [(1)] {\it $R$ is a Hermite ring;}
\vspace{-.5mm}
\item [(2)] {\it Every matrix $\left(
\begin{array}{cc}
a&0\\
b&c
\end{array}
\right)$ with $aR+bR+cR=R$ admits elementary diagonal reduction.}
\end{enumerate}
\vspace{-.5mm} {\it Proof.}\ \ One easily checks that $\left(
\begin{array}{cc}
&1\\
1&
\end{array}
\right)\left(
\begin{array}{cc}
a&0\\
b&c
\end{array}
\right)\left(
\begin{array}{cc}
&1\\
1&
\end{array}
\right)=\left(
\begin{array}{cc}
c&b\\
0&a
\end{array}
\right).$ Therefore the result follows, by [6, Theorem 1.1] and [7, Theorem 2.5].\hfill$\Box$

\vskip4mm \hspace{-1.8em} {\bf Lemma 3.7.}\ \ {\it If $(b+ar)R+cR=R$ with $a,b,c,r\in R$, then $\left(
\begin{array}{cc}
a&b\\
0&c
\end{array}
\right)$ with $aR+bR+cR=R$ admits elementary diagonal reduction.}
\vskip2mm\hspace{-1.8em} {\it Proof.}\ \ Let $A=\left(
\begin{array}{cc}
a&b\\
0&c
\end{array}
\right)$. Since $(b+ar)R+cR=R$, we have $A\left(\begin{array}{cc}
1&r\\
0&1
\end{array}
\right)=\left(\begin{array}{cc}
a&b+ar\\
0&c
\end{array}
\right)$ $=B$. It suffices to prove $B$ admits a diagonal
reduction.
Write $(b+ar)x+cy=1$ for some $x,y\in R$. Then the matrix $\left(\begin{array}{cc}
x&y\\
-c&b+ar
\end{array}
\right)$ is invertible, and we see that  $$\left(\begin{array}{cc}
x&y\\
-c&b+ar
\end{array}
\right)B\left(\begin{array}{cc}
1&\\
-ax&1
\end{array}
\right)\left(\begin{array}{cc}
&1\\
1&
\end{array}
\right)=\left(\begin{array}{cc}
1&0\\
0&-ac
\end{array}
\right),$$ as desired.\hfill$\Box$

\vskip4mm \hspace{-1.8em} {\bf Theorem 3.8.}\ \ {\it If $R$ has feckly adequate range 1, then $R$ is an elementary divisor ring if and only if $R$ is a B$\acute{e}$zout ring.}
\vskip2mm\hspace{-1.8em} {\it Proof.}\ \ $\Longrightarrow $ This is obvious.

$\Longleftarrow$  Step I. Suppose that $aR+bR+cR=R$ with $a,b,c\in R$. Write $ax+by+cz=1$ with $x,y,z\in R$. By hypothesis, there exist $k\in R$ such that
$w:=a+byk+czk\in R$ is feckly adequate. In view of Theorem 3.1, $R/wR$ is feckly zero-adequate. It follows by Theorem 2.5 that $R/wR/J(R/wR)$ is regular, and so it has stable range one. We infer that $R/wR$ has stable range 1. Clearly, $(a+byk+czk)x+by(1-kx)+cz(1-kx)=1$. Thus,
$\overline{by(1-kx)+cz(1-kx)}=\overline{1}$ in $R/wR$. Thus, we can find $h\in R$ such that $\overline{b+cz(1-kx)h}\in U(R/wR)$. It follows that
$\big(b+cz(1-kx)h\big)R+\big(a+byk+czk\big)R=R$. Hence,
$$\big(b+cz(1-kx)h\big)R+\big(a+(b+cz(1-kx)h)yk+czk(1-(1-kx)hy)\big)R=R.$$ Therefore,
$$\big(b+cz(1-kx)h\big)R+\big(a+czk(1-(1-kx)hy)\big)R=R.$$
Thus, $R$ has stable range 2. According to [8, Theorem 2.1.2], $R$ is a Hermite ring.

Step II. Let $A=\left(
\begin{array}{cc}
a'&0\\
b'&c'
\end{array}
\right)$ with $a'R+b'R+c'R=R$. Then there exist $x,y,z\in R$ such that $a'x+b'y+c'z=1$. Since $R$ has feckly adequate range 1,
we can find some $s,t\in R$ such that
$w:=b'+a'xs+c'zt\in R$ is feckly zero-adequate. Hence, $$\left(
\begin{array}{cc}
1&0\\
xs&1
\end{array}
\right)A\left(
\begin{array}{cc}
1&0\\
zt&1
\end{array}
\right)=\left(
\begin{array}{cc}
a'&0\\
w&c'
\end{array}
\right).$$
Since $R$ is a Hermite ring, there exists some $Q=(q_{ij})\in GL_2(R)$ such that $(w,c')Q=(0,c)$ for a $c\in R$.
This implies that $\left(
\begin{array}{cc}
a'&0\\
w&c'
\end{array}
\right)Q=\left(
\begin{array}{cc}
a&b\\
0&c
\end{array}
\right)$.
Clearly, $wR\subseteq cR$. Additionally, we see that $aR+bR+cR=R$. Since $w\in R$ is feckly adequate, $R/wR$ is feckly zero-adequate by
Theorem 3.1. In view of
Theorem 2.5, $(R/wR)/J(R/wR)$ is regular. For any $\overline{\alpha}\in R/wR$, we can find some $\beta \in R$ such that
$\big(\alpha-\alpha\beta\alpha\big)+wR\in J(R/wR)$. It follows that $\big(\alpha-\alpha\beta\alpha\big)+cR\in J(R/cR)$.
Thus, $(R/cR)/J(R/cR)$ is regular; hence, $(R/cR)/J(R/cR)$ has stable range 1. This implies that $R/cR$ has stable range 1.

Clearly,
$\overline{a}(R/cR)+\overline{b}(R/cR)=R/cR$. Then we can
find some $r\in R$ such that $\overline{b+ar}\in R/cR$ is
invertible. Hence, $\overline{(b+ar)d}=\overline{1}$, and then
$(b+ar)d+cp=1$ for some $p\in R$. Therefore $(b+ar)R+cR=R$.
In light of Lemma 3.6, $\left(
\begin{array}{cc}
a&b\\
0&c
\end{array}
\right)$ admits elementary diagonal reduction, and therefore so does $A$, as asserted.
\hfill$\Box$

\vskip4mm \hspace{-1.8em} {\bf Corollary 3.9.}\ \ {\it Let $R$ be a B$\acute{e}$zout ring. If every $a\not\in J(R)$ is feckly adequate,
then $R$ is an elementary divisor ring.}
\vskip2mm\hspace{-1.8em} {\it Proof.}\ \ Suppose that $pR+qR=R$ with $p,q\in R$. If $p\in J(R)$, then $q\in U(R)$. Hence, $p+q\in U(R)$, and so $p+q\in R$ is feckly adequate.
If $p\not\in J(R)$, then $p+q\cdot 0\in R$ is feckly adequate.
Therefore $R$ has feckly adequate range 1, and then we obtain the result, by Theorem 3.8.\hfill$\Box$.

\vskip4mm A B$\acute{e}$zout ring $R$ is called feckly adequate ring provided that every nonzero element in $R$ is feckly adequate. As an immediate consequence of Corollary 3.9, we conclude that every feckly adequate ring is an elementary divisor ring. This generalizes [11, Theorem 7] as well.

\vskip4mm \hspace{-1.8em} {\bf Example 3.10.}\ \ {\it Let $R= \{ a +
bx~|~a\in {\Bbb Z}, b\in {\Bbb Q},x^2=0\}$. Then $R$ is a B$\acute{e}$zout ring in which every $a\not\in J(R)$ is feckly adequate.}
\vskip2mm\hspace{-1.8em} {\it Proof.}\ \ Let $J=(a_1+b_1x)R+(a_2+b_2x)R$. Set $I=\{ \alpha\in {\Bbb Z}~\mid~\alpha+\beta x\in J$ for some $\beta\in {\Bbb Q}\}$. Since ${\Bbb Z}$ is a principal ideal domain, we have some $p,q\in {\Bbb Z}$ such that
$a_1{\Bbb Z}+a_2{\Bbb Z}=p{\Bbb Z}$ and $b_1{\Bbb Z}+b_2{\Bbb Z}=q{\Bbb Z}$. If $I\neq 0$, then $J=pR$. If $I=0$, then $J=(qx)R$. Thus, $R$ is a B$\acute{e}$zout ring. Clearly, $J(R)=x{\Bbb Q}$. Let $f(x)=y+bx\not\in J(R)$, and let $h(x)=z+cx\in R$. Then $y\neq 0$.
Since ${\Bbb Z}$ is a principal ideal domain, it is adequate. Thus, there exist $s,t\in R$ such that
$y=st, (s,z)=1$, and that $(t',z)\neq 1$ for any non-unit divisor $t'$ of $t$.
If $(s,t)\neq 1$, then we have a nonunit $d\in R$ such that $(s,t)=d$. Hence,
$(d,z)\neq 1$, and then $(s,z)\neq 1$, an absurd. Therefore $(s,t)=1$, and so
we can find some $e,d\in R$ such that $se+dt=b.$ One easily checks that
$f(x)=\big(s+dx)\big(t+ex)$. Set $s(x)=s+dx$ and $t(x)=t+ex$. Then $f(x)=s(x)t(x)$.
Clearly, we can find some $k,l\in {\Bbb Z}$ such that $ks+lz=1$. Hence, $1-\big(ks(x)+lg(x)\big)\in J(R)$. Thus,
$ks(x)+lh(x)\in U(R)$. This shows that $\big(s(x),h(x))=1$. If $t'(x)=m+fx$ is a nonunit divisor of $t(x)$, then $m$ is a nonunit divisor of $t$.
By hypothesis, $(m,z)\neq 1$. This implies that $\big(t'(x),h(x)\big)\neq 1$. Thus,
$f(x)\in R$ is adequate, and so $f(x)$ is feckly adequate. In this case, $R$ is not feckly zero-adequate.\hfill$\Box$

\vskip4mm Following Domsha and Vasiunyk, a ring $R$ has adequate range 1 provided that $aR+bR=R$ implies that there exist a $y\in R$ such that $a+by\in R$ is adequate ([4]). For instance, every VNL ring (i.e. for any $a\in R$, either $a$ or $1-a$ is regular) has adequate range 1 ([4, Theorem 11 and Theorem 12]).
We now extend [4, Theorem 14] to B$\acute{e}$zout rings (maybe with zero divisors)

\vskip4mm \hspace{-1.8em} {\bf Corollary 3.11.}\ \ {\it If $R$ has adequate range 1, then
$R$ is an elementary divisor ring if and only if $R$ is a B$\acute{e}$zout ring.}
\vskip2mm\hspace{-1.8em} {\it Proof.}\ \ As every adequate element in a ring is feckly adequate, if $R$ has adequate range 1, then it has feckly adequate range 1. Therefore we complete the proof, by Theorem 3.8.\hfill$\Box$.

\vskip20mm \bc{\bf References}\ec \vskip4mm {\small \re{1} D.D. Anderson and V.P. Camillo,
Commutative rings whose elements are a sum of a unit and
idempotent, {\it Comm. Algebra}, {\bf 30}(2002), 3327--3336.

\re{2} H.
Chen, {\it Rings Related Stable Range Conditions}, Series in
Algebra 11, World Scientific, Hackensack, NJ, 2011.

\re{3} H. Chen; H. Kose and Y. Kurtulmaz, On feckly clean rings, {\it J. Algebra Appl.}, 11/2014;  DOI: 10.1142/S0219498815500462.

\re{4} O.V. Domsha and I.S. Vasiunyk, Combing local and adequate rings, {\it Book of abstracts of the International Algebraic Conference}, Taras Shevchenko National University of Kyiv, Kyiv, Ukraine, 2014.

\re{5} O. Pihura, Commutative Bezout rings in which 0 is adequate is a semi regular, preprint, 2015.1.

\re{6} W.W. McGovern, B\'{e}zout rings with almost stable range
$1$, {\it J. Pure Appl. Algebra}, {\bf 212}(2008), 340--348.

\re{7} M. Roitman, The Kaplansky condition and rings of almost
stable range $1$, {\it Proc. Amer. Math. Soc.}, {\bf 141}(2013),
3013--3018.

\re{8} B.V. Zabavsky, Diagonal Reduction of Matrices over Rings,
Mathematical Studies Monograph Series, Vol. XVI, VNTL Publisher, 2012.

\re{9} B.V. Zabavsky, Questions related to the K-theoretical
aspect of B$\acute{e}$zout rings with various stable range conditions,
{\it Math. Stud.}, {\bf 42}(2014), 89-103.

\re{10} B.V. Zabavsky and S.I. Bilavska,
Every zero adequate ring is an exchange ring, {\it J. Math. Sci.},
{\bf 187}(2012), 153--156.

\re{11} B.V. Zabavsky and M.Y. Komarnyts'kyi, Cohn-type theorem
for adequacy and elementary divisor rings, {\it J. Math. Sci.},
{\bf 167}(2010), 107--111.

\end{document}